\documentclass[12pt,english,twoside]{article}
\usepackage[english]{babel}
\usepackage{amsmath,amsfonts}
\usepackage{mathtools}
\usepackage[cp850]{inputenc}
\allowdisplaybreaks
\newcounter{theorem}[section]
\renewcommand{\thetheorem}{\thesection.\arabic{theorem}}
\newenvironment{lemma}[      1]{\refstepcounter{theorem} %
\bf \thetheorem\ Lemma#1.        \it}{}
\newenvironment{theorem}[    1]{\refstepcounter{theorem} %
\bf \thetheorem\ Theorem#1.      \it}{}

\newenvironment{corollary}[  1]{\refstepcounter{theorem} %
\bf \thetheorem\ Corollary#1.    \it}{}
\newenvironment{example}[    1]{\refstepcounter{theorem} %
\bf \thetheorem\ Example#1.      \rm}{\par}
\newenvironment{examples}[   1]{\refstepcounter{theorem} %
\bf \thetheorem\ Examples#1.     \rm}{\par}

\newenvironment{proof}{ %
\it              Proof.          \rm}{\hfill $ \Box $}

\newcounter{abc}[theorem]
\newenvironment{abclist}{\begin{list}{%
\rm (\alph{abc})   \hfill        \it}{\usecounter{abc} %
\topsep0mm \partopsep0mm \parsep0mm \itemsep0mm %
\leftmargin2em \labelwidth2em \labelsep0em}}{\end{list}}

\newcounter{one}[theorem]
\newenvironment{onelist}{\begin{list}{%
\rm (\arabic{one}) \hfill           }{\usecounter{one} %
\topsep0mm \partopsep0mm \parsep0mm \itemsep0mm %
\leftmargin2em \labelwidth2em \labelsep0em}}{\end{list}}

\newenvironment{hylist}{\begin{list}{%
--               \hfill            }{%
\topsep0mm \partopsep0mm \parsep0mm \itemsep0mm %
\leftmargin2em \labelwidth2em \labelsep0em}}{\end{list}}

\newcounter{rom}

\newcommand{\N}{{\mathbb N}}
\newcommand{\R}{{\mathbb R}}

\newcommand{\BB}{{\cal B}}

\newcommand{\FF}{{\cal F}}
\newcommand{\LL}{{\cal L}}

\newcommand{\one}{{\prime}}

\newcommand{\leb}{{\lambda}}

\newcommand{\VaR}{{\text{\rm VaR}}}
\newcommand{\ES}{{\text{\rm ES}}}

\newcommand{\Acerbi}{\text{\rm Acerbi}}
\newcommand{\Pichler}{\text{\rm Pichler}}

\newcommand{\leftmatrix}[1]{\left(\begin{array}{#1}}
\newcommand{\rightmatrix}{\end{array}\right)}

\begin{document}

\title{\Large\bf On Quantile Risk Measures and Their Domain}
\author{Sebastian Fuchs\footnote{Faculty of Economics and Management, Free University of Bozen--Bolzano, 39100 Bozen, Italy}, 
Ruben Schlotter\footnote{Fakult{\"a}t f{\"u}r Mathematik, Technische Universit{\"a}t Chemnitz, 09126 Chemnitz, Germany},  
and 
Klaus D. Schmidt\footnote{Fachrichtung Mathematik, Technische Universit{\"a}t Dresden, 01062 Dresden, Germany}}
\date{}
\maketitle

%\title{\bf Bivariate Copulas: Transformations, Asymmetry and Measures of Concordance}
%\author{Sebastian Fuchs}
%\author{Klaus D. Schmidt\corref{cor1}}
%\ead{klaus.d.schmidt@tu-dresden.de}
%\address{Lehrstuhl f\"ur Versicherungsmathematik, Technische Universit\"at Dresden, \\ 01062 Dresden, Germany.}

\begin{abstract}
\noindent
In the present paper we study quantile risk measures and their domain. 
Our starting point is that, 
for a probability measure   $ Q $   on the open unit interval and a wide class   $ \LL_Q $   of random variables, 
we define the quantile risk measure   $ \varrho_Q $   as the map 
which integrates the quantile function of a random variable in   $ \LL_Q $   with respect to   $ Q $. 
The definition of   $ \LL_Q $   ensures that   $ \varrho_Q $   
cannot attain the value   $ +\infty $   and 
cannot be extended beyond   $ \LL_Q $   without losing this property. 
The notion of a quantile risk measure is a natural generalization of that of a spectral risk measure 
and provides another view at the distortion risk measures generated by a distribution function on the unit interval. 
In this general setting, 
we prove several results on quantile or spectral risk measures and their domain 
with special consideration of the expected shortfall. 
\end{abstract}

%\begin{keyword}
%Copulas \sep transformations \sep involutions \sep asymmetry \sep invariance \sep order \sep measures of concordance 
%\end{keyword}
%\subclass{Primary 65H10, 47J10; Secondary 91B30.}
%\maketitle

%%%%%%%%%%%%%%%%%%%%%%%%%%%%%%%%%%%%%%%%%%%%%%%%%%%%%%%%%%%%%%%%%%%%%%%%%%%%%%%%%%%%%%%%%%%%%%%%%%%

\section{Introduction}
\label{introduction}

In the present paper we study quantile risk measures and their domain. 
Our starting point is that, 
for a probability measure   $ Q $   on the open unit interval and a wide class   $ \LL_Q $   of random variables, 
we define the quantile risk measure   $ \varrho_Q $   as the map 
which integrates the quantile function of a random variable in   $ \LL_Q $   with respect to   $ Q $. 
The definition of   $ \LL_Q $   ensures that   $ \varrho_Q $   
cannot attain the value   $ +\infty $   and 
cannot be extended beyond   $ \LL_Q $   without losing this property. 
The notion of a quantile risk measure is a natural generalization of that of a spectral risk measure 
and provides another view at the distortion risk measures generated by a distribution function on the unit interval. 

\bigskip
Quantile risk measures are thus 
mixtures of the values at risk at different levels and hence 
mixtures of a parametric family of risk measures. 
Such mixtures have already been considered by 
Acerbi [2002] who, 
however, 
spent little attention to the domain on which a given risk measure can be defined; 
he argued that 
\emph{in a real--world risk management application the integral} (defining a risk measure) 
\emph{will always be well--defined and finite}. 
Nevertheless, 
%it should be pointed out that 
Acerbi [2002] 
%he 
proposed a maximal class of random variables 
on which a given spectral risk measure is well--defined and finite. 
In the case of a spectral risk measure, 
the domain of a quantile risk measure proposed in the present paper contains the class proposed by 
Acerbi [2002] and turns out to be a convex cone, 
which is of interest with regard to subadditivity of the risk measure. 

\bigskip
In this paper we review and partly extend known results on quantile risk measures, 
with particular attention to spectral risk measures and, 
in particular, 
expected shortfall, 
and with emphasis on their maximal domain mentioned before.  
We deliberately adopt arguments from the literature, 
with appropriate modifications if necessary, 
but some of our proofs and results are new. 

%\bigskip
%Sloppiness with regard to the domain is widespread in the literature on risk measures: 
%Certain authors do not define the domain at all, 
%others assume that all random variables in the domain are bounded and/or positive, 
%and some authors even assume that the risk measure under consideration is well--defined for every random variable 
%(which is true for value at risk but is false for the expectation). 
%Besides the paper by 
%Acerbi [2002], 
%a notable exception is the recent paper by 
%Pichler [2013] who proposed a vector lattice of random variables 
%on which a given spectral risk measure is well--defined and finite. 
%
%\bigskip
%In view of these observations, 
%we found that the time has come to study risk measures on a maximal set of random variables 
%for which they are well--defined and cannot attain the value   $ +\infty $. 
%In the present paper, 
%we carry out this programme for quantile risk measures, 
%which essentially correspond to distortion risk measures, 
%for the special case of spectral risk measures and, 
%in particular, 
%for expected shortfall. 
%In doing so, 
%we extend results which are known for certain subsets of the domain of a quantile risk measure, 
%and in the proofs we deliberately adopt arguments from the literature with appropriate modifications. 

\bigskip
This paper is organized as follows: 
We first fix some notation and recall some elementary properties of the quantile function 
and the basic examples of distortion functions 
(Section \ref{preliminaries}). 
We then introduce quantile risk measures and provide several alternative representations of quantile risk measures and their domain, 
as well as 
%a condition for a quantile risk measure to be finite and 
conditions under which certain quantile risk measures can be compared 
(Section \ref{quantile}). 
In the next step, 
we introduce spectral risk measures and characterize spectral risk measures within the class of all quantile risk measures
(Section \ref{spectral}). 
We then present a short proof of the subadditivity of expected shortfall and use this result to show 
that a quantile risk measures is subadditive if and only if it is 
spectral 
%a spectral risk measure   
(Section~\ref{subadditive}). 
We conclude with a comparison of 
the domain of a quantile risk measure with the classes of random variables proposed by 
Acerbi [2002] and 
Pichler [2013] in the spectral case 
(Section~\ref{domain}).

%%%%%%%%%%%%%%%%%%%%%%%%%%%%%%%%%%%%%%%%%%%%%%%%%%%%%%%%%%%%%%%%%%%%%%%%%%%%%%%%%%%%%%%%%%%%%%%%%%%%%%%%

\section{Preliminaries}
\label{preliminaries}

We use the terms 
\emph{positive} and 
\emph{increasing} in the weak sense which admits equality in the inequalities defining these terms. 
For   $ B\subseteq\R $, 
we denote by   $ \chi_B $   the indicator function of   $ B $   
(such that   
$ \chi_B(x)=1 $   if   $ x \in B $   and 
$ \chi_B(x)=0 $   if   $ x \notin B $). 
Also, 
we denote 
\begin{hylist}
\item   by   $ \BB(\R) $   the $\sigma$--field of all Borel sets of   $ \R $, 
\item   by   $ \BB((0,1)) $   the $\sigma$--field of all Borel sets of   $ (0,1) $,   and 
\item   by   $ \leb $   the Lebesgue measure on   $ \BB(\R) $   or its restriction to   $ \BB((0,1)) $. 
\end{hylist}
By the correspondence theorem, 
there exists a bijection between the distribution functions on   $ \R $      and the probability measures on   $ \BB(\R) $   such that 
the probability measure   $ Q^G $   corresponding to the distribution function   $ G $   satisfies   $ Q^G[(x,y]] = G(y) - G(x) $   
for all  $ x,y\in\R $   such that   $ x \leq y $. 
Correspondingly, 
there exists a bijection between the distribution functions on   $ (0,1) $   and the probability measures on   $ \BB((0,1)) $. 

\bigskip
Throughout this paper, 
we consider a fixed probability space   $ (\Omega,\FF,P) $   and random variables   $ (\Omega,\FF)\to(\R,\BB(\R)) $   and we denote 
\begin{hylist}
\item   by   $ \LL^0 $   the vector lattice of all random variables, 
\item   by   $ \LL^1 $   the vector lattice of all integrable random variables, 
\item   by   $ \LL^2 $   the vector lattice of all square integrable random variables, and 
\item   by   $ \LL^\infty $   the vector lattice of all almost surely bounded random variables. 
\end{hylist}
Then we have   $ \LL^\infty\subseteq\LL^2\subseteq\LL^1\subseteq\LL^0 $. 
For a random variable   $ X\in\LL^0 $   we denote 
by   $ F_X $   its distribution function   $ \R\to[0,1] $   given by 
\begin{eqnarray*}
        F_X(x)           
& := &  P[\{X \leq x\}]  
\end{eqnarray*}
and 
by   $ F^\leftarrow_X $   its (lower) quantile function    $ (0,1)\to\R $   given by 
\begin{eqnarray*}
        F^\leftarrow_X(u)                                
& := &  \inf\Bigl\{ x\in\R \Bigm| F_X(x) \geq u \Bigr\}  
\end{eqnarray*}
For   $ u\in(0,1) $   and   $ x\in\R $, 
the quantile function satisfies   $ F^\leftarrow_X(u) \leq x $   if and only if   $ u \leq F_X(x) $. 
Moreover, 
the quantile function is increasing and has the following properties: 

\bigskip
\begin{lemma}{}
\label{preliminaries-l}
Consider   $ X,Y\in\LL^0 $. 
Then: 
\begin{onelist}
\item   If   $ X \leq Y $,   
then   $ F^\leftarrow_X \leq F^\leftarrow_Y $. 
\item   If   $ a\in\R_+ $, 
then   $ F^\leftarrow_{aX} = a\,F^\leftarrow_X $. 
\item   If   $ c\in\R $, 
then   $ F^\leftarrow_{X+c} = F^\leftarrow_X + c $. 
\item   If   $ X $   and   $ Y $   are comonotone, 
then   $ F^\leftarrow_{X+Y} = F^\leftarrow_X + F^\leftarrow_Y $. 
\item   $ F^\leftarrow_{X^+} = (F^\leftarrow_X)^+ $. 
\end{onelist}
\end{lemma}

\bigskip
A function   $ D : [0,1]\to[0,1] $   is said to be a 
\emph{distortion function} if it is increasing and continuous from the right 
and satisfies   $ D(0) = 0 $   and   $ \sup_{u\in(0,1)} D(u) = 1 $   (and hence   $ D(1) = 1 $). 
The restriction of a distortion function   $ D $   to   $ (0,1) $   is a distribution function on   $ (0,1) $   and, 
for simplicity,  
the probability measure corresponding to the restriction of   $ D $   to   $ (0,1) $   will be referred to as 
the probability measure corresponding to   $ D $. 

\bigskip
\begin{examples}{}
\label{preliminaries-e}
The terms attached to the following examples are the names of the risk measures resulting from the respective distortion functions. 
\begin{onelist}
\item   {\bf Expectation:} 
The function   $ D^E : [0,1]\to[0,1] $   given by 
\begin{eqnarray*}
        D^E(u)  
& := &  u       
\end{eqnarray*}
is a distortion function. 
\item   {\bf Value at Risk:} 
For   $ \alpha\in(0,1) $, 
the function   $ D^{\VaR_\alpha} : [0,1]\to[0,1] $   given by 
\begin{eqnarray*}
        D^{\VaR_\alpha}(u)    
& := &  \chi_{[\alpha,1]}(u)  
\end{eqnarray*}
is a distortion function. 
\item   {\bf Expected Shortfall:} 
For   $ \alpha\in[0,1) $, 
the function   $ D^{\ES_\alpha} : [0,1]\to[0,1] $   given by 
\begin{eqnarray*}
        D^{\ES_\alpha}(u)                                
& := &  \frac{u-\alpha}{1-\alpha}\,\chi_{[\alpha,1]}(u)  
\end{eqnarray*}
is a distortion function; 
in particular, 
$ D^{\ES_0} = D^E $. 
\item   {\bf Expected Shortfall of Higher Order:} 
For   $ n\in\N $   and   $ \alpha\in[0,1) $, 
the function   $ D^{\ES_{n,\alpha}}(u) : [0,1]\to[0,1] $   given by 
\begin{eqnarray*}
        D^{\ES_{n,\alpha}}(u)                                              
& := &  \biggl( \frac{u-\alpha}{1-\alpha} \biggr)^n\,\chi_{[\alpha,1]}(u)  
\end{eqnarray*}
is a distortion function; 
in particular, 
$ D^{\ES_{1,\alpha}} = D^{\ES_\alpha} $. 
\end{onelist}
The distortion functions   $ D^{\ES_{n,\alpha}} $, 
and in particular   $ D^{\ES_\alpha} $   and   $ D^E$,   
are convex whereas 
the distortion functions   $ D^{\VaR_\alpha} $   are not convex. 
\end{examples}

\bigskip
Throughout this paper, 
we consider pairs   $ (D,Q) $   consisting of a distortion function   $ D : [0,1]\to[0,1] $   
and the probability measure   $ Q : \BB((0,1))\to[0,1] $   corresponding to   $ D $, 
and we use identical sub-- or superscripts for both, 
$ D $   and   $ Q $, 
in the case of a particular choice of   $ D $   or   $ Q $.

%%%%%%%%%%%%%%%%%%%%%%%%%%%%%%%%%%%%%%%%%%%%%%%%%%%%%%%%%%%%%%%%%%%%%%%%%%%%%%%%%%%%%%%%%%%%%%%%%%%%%%%%

\section{Quantile Risk Measures}
\label{quantile}

Define 
\begin{eqnarray*}
        \LL_Q                                                                                   
& := &  \biggl\{ X\in\LL^0 \biggm| \int_{(0,1)} (F^\leftarrow_X(u))^+\,dQ(u) < \infty \biggr\}  
\end{eqnarray*}
Then we have   $ \LL^\infty\subseteq\LL_Q $   and the map   $ \varrho_Q : \LL_Q\to[-\infty,\infty) $   given by 
\begin{eqnarray*}
        \varrho_Q[X]                           
& := &  \int_{(0,1)} F^\leftarrow_X(u)\,dQ(u)  
\end{eqnarray*}
is said to be a 
\emph{quantile risk measure}. 

\bigskip
For every   $ X\in\LL^0 $,   
we have   $ X\in\LL_Q $   if and only if   $ X^+\in\LL_Q $, 
by Lemma 
\ref{preliminaries-l}. 
This implies that, 
for every   $ Z\in\LL^0  $   satisfying   $ Z \leq X $   for some   $ X\in\LL_Q $, 
we have   $ Z\in\LL_Q $. 
Lemma 
\ref{preliminaries-l} also yields the following properties of a quantile risk measure: 

\bigskip
\begin{lemma}{}
\label{quantile-l}
Consider   $ X,Y\in\LL_Q $. 
Then: 
\begin{onelist}
\item   If   $ X \leq Y $,   
then   $ \varrho_Q[X] \leq \varrho_Q[Y] $. 
\item   If   $ a\in\R_+ $, 
then   $ aX\in\LL_Q $   and   $ \varrho_Q[aX] = a\,\varrho_Q[X] $. 
\item   If   $ c\in\R $, 
then   $ X\!+\!c\in\LL_Q $   and   $ \varrho_Q[X\!+\!c] = \varrho_Q[X] + c $. 
\item   If   $ X $   and   $ Y $   are comonotone, 
then   $ X+Y\in\LL_Q $   and   $ \varrho_Q[X+Y] = \varrho_Q[X] + \varrho_Q[Y] $. 
\end{onelist}
\end{lemma}

\bigskip
The quantile risk measure   $ \varrho_Q $   is said to be 
%\begin{hylist}
%\item   \emph{additive}    if   $ \varrho_Q[X\!+\!Y]  =   \varrho_Q[X] + \varrho_Q[Y] $   holds for all   $ X,Y\in\LL_Q $   such that   $ X+Y\in\LL_Q $. 
%\item   \emph{subadditive} if   $ \varrho_Q[X\!+\!Y] \leq \varrho_Q[X] + \varrho_Q[Y] $   holds for all   $ X,Y\in\LL_Q $   such that   $ X+Y\in\LL_Q $. 
%\end{hylist}
\emph{subadditive} if 
$ \varrho_Q[X+Y] \leq \varrho_Q[X] + \varrho_Q[Y] $   
%\begin{eqnarray*}
%        \varrho_Q[X\!+\!Y]           
%&\leq&  \varrho_Q[X] + \varrho_Q[Y]  
%\end{eqnarray*}
holds for all   $ X,Y\in\LL_Q $   such that   $ X+Y\in\LL_Q $. 
We shall show that   $ \varrho_Q $   is subadditive if and only if   $ D $   is convex, 
and that in this case   $ \LL_Q $   is a convex cone; 
see Theorem 
\ref{subadditive-char-t} below. 

\bigskip
To obtain alternative representations of a quantile risk measure and its domain we need the following Lemma: 

\bigskip
\begin{lemma}{}
\label{quantile-char-l}
The identities 
$$  \int_{(0,1)} (F_X^\leftarrow(u))^+\,dQ(u)  =  \int_\R x^+\,dQ^{D \circ F_X}(x)  =  \int_{(0, \infty)} \Bigl( 1- (D \circ F_X)(x) \Bigr)\,d\leb(x)  $$  
and 
$$  \int_{(0,1)} (F_X^\leftarrow(u))^-\,dQ(u)  =  \int_\R x^-\,dQ^{D \circ F_X}(x)  =  \int_{(-\infty,0)}           (D \circ F_X)(x)       \,d\leb(x)  $$  
hold for every   $ X\in\LL^0 $. 
\end{lemma}

\bigskip
\begin{proof}
For every   $ x\in\R $   we have 
\begin{eqnarray*}
        Q^{D \circ F_X}[(-\infty,x]]                                              
&  = &  (D \circ F_X)(x)  	                                                      \\[1ex]
&  = &  Q[(0,F_X(x)]\cap(0,1)]                                                    \\[1ex]
&  = &  Q\Bigl[ \Bigl\{ u\in(0,1) \Bigm| u \leq F_X(x)            \Bigr\} \Bigr]  \\
&  = &  Q\Bigl[ \Bigl\{ u\in(0,1) \Bigm| F^\leftarrow_X(u) \leq x \Bigr\} \Bigr]  \\*[1ex]
&  = &  Q_{F^\leftarrow_X}[(-\infty,x]]
\end{eqnarray*}
and hence   $ Q^{D \circ F_X} = Q_{F^\leftarrow_X} $. 
Now the substitution rule yields 
\begin{eqnarray*}
        \int_{(0,1)} (F^\leftarrow_X(u))^+\,dQ(u)                   
&  = &  \int_\R      x^+                  \,dQ_{F^\leftarrow_X}(x)  \\*
&  = &  \int_\R      x^+                  \,dQ^{D \circ F_X}   (x)  
\end{eqnarray*}
and 
\begin{eqnarray*}
        \int_{(0,1)} (F^\leftarrow_X(u))^-\,dQ(u)                   
&  = &  \int_\R      x^-                  \,dQ_{F^\leftarrow_X}(x)  \\*
&  = &  \int_\R      x^-                  \,dQ^{D \circ F_X}   (x)  
\end{eqnarray*}
Moreover, 
we have 
\begin{eqnarray*}
        \int_\R           x^+                                             \,dQ^{D \circ F_X}(x)  
&  = &  \int_{(0,\infty)} x                                               \,dQ^{D \circ F_X}(x)  \\
&  = &  \int_{(0,\infty)} \int_{(0,\infty)} \chi_{(0,     x)}(y)\,d\leb(y)\,dQ^{D \circ F_X}(x)  \\
&  = &  \int_{(0,\infty)} \int_{(0,\infty)} \chi_{(y,\infty)}(x)\,dQ^{D \circ F_X}(x)\,d\leb(y)  \\*
&  = &  \int_{(0,\infty)} \Bigl( 1 - D \circ F_X(y)\Bigr)                            \,d\leb(y)  
\end{eqnarray*}
and
\begin{eqnarray*}
        \int_\R            x^-                                               \,dQ^{D \circ F_X}(x)  
&  = &  \int_{(-\infty,0)} (-x)                                              \,dQ^{D \circ F_X}(x)  \\
&  = &  \int_{(-\infty,0)} \int_{(-\infty,0)} \chi_{[x      ,0)}(y)\,d\leb(y)\,dQ^{D \circ F_X}(x)  \\
&  = &  \int_{(-\infty,0)} \int_{(-\infty,0)} \chi_{(-\infty,y]}(x)\,dQ^{D \circ F_X}(x)\,d\leb(y)  \\*
&  = &  \int_{(-\infty,0)} (D \circ F_X)(y)                                             \,d\leb(y)  
\end{eqnarray*}
The assertion follows. 
\end{proof}

\bigskip
The following result is immediate from Lemma 
\ref{quantile-char-l}: 

\bigskip
\begin{theorem}{}
\label{quantile-char-t}
The domain of   $ \varrho_Q $   satisfies 
\begin{eqnarray*}
        \LL_Q                                                                                                                   
&  = &  \biggl\{ X \in \LL^0 \biggm| \int_\R x^+\,dQ^{D \circ F_X}(x)                                        < \infty \biggr\}  \\*
&  = &  \biggl\{ X \in \mathcal{L}^{0} \biggm| \int_{(0,\infty)} \Bigl( 1- (D \circ F_X)(x) \Bigr)\,d\leb(x) < \infty \biggr\}  
\end{eqnarray*}
and the identities 
\begin{eqnarray*}
        \varrho_Q[X]                                                                                                    
&  = &  \int_\R x\,dQ^{D \circ F_X}(x)                                                                                  \\*
&  = &  \int_{(0,\infty)} \Bigl( 1 - (D \circ F_X)(x) \Bigr)\,d\leb(x) - \int_{(-\infty,0)} (D \circ F_X)(x)\,d\leb(x)  
\end{eqnarray*}
hold for every   $ X \in \LL_Q $. 
\end{theorem}

\bigskip
Because of the previous result, 
the quantile risk measure generated by the probability measure   $ Q $   corresponds to 
the 
\emph{distortion risk measure} generated by the distortion function   $ D $; 
the latter is also known as 
\emph{Wang's premium principle}.  

\bigskip
\begin{examples}{}
\label{quantile-e}
\begin{onelist}
\item   {\bf Expectation:} 
The distortion function   $ D^E $   satisfies   $ D^E \circ F_X = F_X $. 
Because of Theorem 
\ref{quantile-char-t} this yields 
\begin{eqnarray*}
        \LL_{Q^E}                                         
&  = &  \Bigl\{ X\in\LL^0 \Bigm| E[X^+] < \infty \Bigr\}  
\end{eqnarray*}
and 
\begin{eqnarray*}
        \varrho_{Q^E}[X]  
&  = &  E[X]              
\end{eqnarray*}
for every   $ X\in\LL_{Q^E} $. 
\item   {\bf Value at Risk:} 
For   $ \alpha\in(0,1) $, 
the probability measure   $ Q^{\VaR_\alpha} $   corresponding to   $ D^{\VaR_\alpha} $   is the Dirac measure at   $ \alpha $. 
This yields 
\begin{eqnarray*}
        \LL_{Q^{\VaR_\alpha}}  
&  = &  \LL^0                  
\end{eqnarray*}
(such that   $ \LL_{Q^{\VaR_\alpha}} $   does not depend on   $ \alpha $) and 
\begin{eqnarray*}
        \varrho_{Q^{\VaR_\alpha}}[X]  
&  = &  F^\leftarrow_X(\alpha)        
\end{eqnarray*}
for every   $ X\in\LL_{Q^{\VaR_\alpha}} $; 
in particular, 
$ \varrho_{Q^{\VaR_\alpha}} $   is finite. 
The quantile risk measure   $ \varrho_{Q^{\VaR_\alpha}} $   is called 
\emph{value at risk at level   $ \alpha $} and is usually denoted by   $ \VaR_\alpha $. 
\item   {\bf Expected Shortfall:} 
For   $ \alpha\in[0,1) $, 
the probability measure   $ Q^{\ES_\alpha} $   corresponding to   $ D^{\ES_\alpha} $   satisfies 
\begin{eqnarray*}
        Q^{\ES_\alpha}                                           
&  = &  \int \frac{1}{1-\alpha}\,\chi_{(\alpha,1)}(u)\,d\leb(u)  
\end{eqnarray*}
Since   $ F^\leftarrow_X $   is increasing and   $ F^\leftarrow_X(\alpha) $   is finite for   $ \alpha\in(0,1) $, 
this yields, 
because of (1),  
\begin{eqnarray*}
        \LL_{Q^{\ES_\alpha}}                                                                                  
&  = &  \biggl\{ X \in \LL^0 \biggm| \int_{(\alpha,1)} (F^\leftarrow_X(u))^+ \,d\lambda(u) < \infty \biggr\}  \\
&  = &  \biggl\{ X \in \LL^0 \biggm| \int_{(     0,1)} (F^\leftarrow_X(u))^+ \,d\lambda(u) < \infty \biggr\}  \\[.5ex]
&  = &  \Bigl\{ X\in\LL^0 \Bigm| E[X^+] < \infty \Bigr\}                                                      \\*[1ex]
&  = &  \LL^E                                                                                                 
\end{eqnarray*}
(such that   $ \LL_{Q^{\ES_\alpha}} $   does not depend on   $ \alpha $) and 
\begin{eqnarray*}
        \varrho_{Q^{\ES_\alpha}}[X]                                                         
&  = &  \int_{(0,1)} F^\leftarrow_X(u)\,\frac{1}{1-\alpha}\,\chi_{(\alpha,1)}(u)\,d\leb(u)  
\end{eqnarray*}
for every   $ X\in\LL_{Q^{\ES_\alpha}} $. 
In particular, 
$ \varrho_{Q^{\ES_0}} = \varrho_{Q^E} $   
and   $ \varrho_{Q^{\ES_\alpha}} $   is finite for every   $ \alpha\in(0,1) $. 
The quantile risk measure   $ \varrho_{Q^{\ES_\alpha}} $   is called 
\emph{expected shortfall at level   $ \alpha $} and is usually denoted by   $ \ES_\alpha $. 
\item   {\bf Expected Shortfall of Higher Order:} 
For   $ n\in\N $   and   $ \alpha\in[0,1) $, 
the probability measure   $ Q^{\ES_{n,\alpha}} $   corresponding to   $ D^{\ES_{n,\alpha}} $   satisfies 
\begin{eqnarray*}
        Q^{\ES_{n,\alpha}}                                                                                       
&  = &  \int \frac{n}{1-\alpha} \biggl( \frac{u-\alpha}{1-\alpha} \biggr)^{n-1}\,\chi_{(\alpha,1)}(u)\,d\leb(u)  
\end{eqnarray*}
This yields 
\begin{eqnarray*}
        \LL_{Q^{\ES_{n,\alpha}}}  
&  = &  \LL_{Q^E}                 
\end{eqnarray*}
(such that   $ \LL_{Q^{\ES_{n,\alpha}}} $   does not depend on   $ n $   or   $ \alpha $) and 
\begin{eqnarray*}
        \varrho_{Q^{\ES_{n,\alpha}}}[X]                                                                                                     
&  = &  \int_{(0,1)} F^\leftarrow_X(u)\,\frac{n}{1-\alpha} \biggl( \frac{u-\alpha}{1-\alpha} \biggr)^{n-1}\,\chi_{(\alpha,1)}(u)\,d\leb(u)  
\end{eqnarray*}
for every   $ X\in\LL_{Q^{\ES_{n,\alpha}}} $. 
In particular, 
$ \varrho_{Q^{\ES_{1,\alpha}}} = \varrho_{Q^{\ES_\alpha}} $   
and   $ \varrho_{Q^{\ES_{n,\alpha}}} $   is finite for every   $ \alpha\in(0,1) $. 
The quantile risk measure   $ \varrho_{Q^{\ES_{n,\alpha}}} $   is called 
\emph{expected shortfall of order   $ n $   at level   $ \alpha $}. 
\end{onelist}
\end{examples}

\bigskip
Lemma 
\ref{quantile-char-l} and Theorem 
\ref{quantile-char-t} have several applications. 
For example, 
they provide a condition on   $ D $   under which   $ \varrho_Q $   is finite: 

\bigskip
\begin{corollary}{}
\label{quantile-finite-c}
Assume that there exists some   $ \delta\in(0,1) $   such that   $ D(u)=0 $   holds for every   $ u\in(0,\delta) $. 
Then 
\begin{eqnarray*}
        \LL_Q                                                                                       
&  = &  \biggl\{ X\in\LL^0 \biggm| \int_{(0,1)} |F^\leftarrow_X(u)|\,dQ(u)       < \infty \biggr\}  \\
&  = &  \biggl\{ X\in\LL^0 \biggm| \int_\R |x|\,dQ^{D \circ F_X}(x)              < \infty \biggr\}  \\*
&  = &  \biggl\{ X\in\LL^0 \biggm| \int_{(0,\infty)} \Bigl( 1- (D \circ F_X)(x) \Bigr)\,d\leb(x)    
                                 + \int_{(-\infty,0)} (D \circ F_X)(x)\,d\leb(x) < \infty \biggr\}  
\end{eqnarray*}
and   $ \varrho_Q $   is finite. 
\end{corollary}

\bigskip
\begin{proof}
For every   $ X\in\LL^0 $, 
the assumption yields 
\begin{eqnarray*}
        \int_{(0,1)} (F^\leftarrow_X(u))^-                                       \,dQ(u)     
&  = &  \int_{(-\infty,0)} (D \circ F_X)(x)                                      \,d\leb(x)  \\
&  = &  \int_{(-\infty,0)} (D \circ F_X)(x)\,\chi_{[\delta,1)}(F_X(x))           \,d\leb(x)  \\
&  = &  \int_{(-\infty,0)} (D \circ F_X)(x)\,\chi_{[F_X^\leftarrow(\delta),0)}(x)\,d\leb(x)  \\*[1ex]
&\leq&  (D \circ F_X)(0) \int_{(-\infty,0)}  \chi_{[F_X^\leftarrow(\delta),0)}(x)\,d\leb(x)  
\end{eqnarray*}
Since   $ F_X^\leftarrow(\delta) $   is finite, 
this proves the assertion. 
\end{proof}

\bigskip
Theorem 
\ref{quantile-char-t} also provides a condition for the comparison of the domains of quantile risk measures: 

\bigskip
\begin{corollary}{}
\label{quantile-order-1-c}
Assume that there exists some   $ \delta\in(0,1) $   such that   $ D_1(u) \leq D_2(u) $   holds for every   $ u\in[\delta,1) $. 
Then   $ \LL_{Q_1} \subseteq \LL_{Q_2} $.
\end{corollary}

\bigskip
\begin{proof}{}
For every   $ X\in\LL^0 $, 
we have 
\begin{eqnarray*}
        \int_{(0,\infty)} \Bigl( 1 - (D_2 \circ F_X)(x) \Bigr)                                           \,d\leb(x)  
&  = &  \int_{(0,\infty)} \Bigl( 1 - (D_2 \circ F_X)(x) \Bigr)\,\chi_{(0, F_X^\leftarrow(\delta)    )}(x)\,d\leb(x)  \\*
&&   +\,\int_{(0,\infty)} \Bigl( 1 - (D_2 \circ F_X)(x) \Bigr)\,\chi_{[F_X^\leftarrow(\delta),\infty)}(x)\,d\leb(x)  \\
&\leq&  \int_{(0,\infty)}                                       \chi_{(0, F_X^\leftarrow(\delta)    )}(x)\,d\leb(x)  \\*
&&   +\,\int_{(0,\infty)} \Bigl( 1 - (D_1 \circ F_X)(x) \Bigr)                                           \,d\leb(x)  
\end{eqnarray*}
Since   $ F_X^\leftarrow(\delta) $   is finite, 
Theorem 
\ref{quantile-char-t} yields   $ \LL_{Q_1}\subseteq\LL_{Q_2} $. 
\end{proof}

\bigskip
\begin{corollary}{}
\label{quantile-order-ES-c} 
Assume that there exist some   $ n\in\N $   and   $ \alpha,\delta\in(0,1) $   such that 
$$  D^{\ES_{n,\alpha}}(u)  \leq  D(u)  \leq  D^E(u)  $$
holds for every   $ u\in[\delta,1) $.
Then   $ \LL_Q = \LL_{Q^E} $.
\end{corollary}

\bigskip
\begin{proof}
Because of Corollary 
\ref{quantile-order-1-c}, 
we have   $ \LL_{Q^{\ES_{n,\alpha}}}\subseteq\LL_Q\subseteq\LL_{Q^E} $. 
Now the assertion follows from   $ \LL_{Q^{\ES_{n,\alpha}}}=\LL_{Q^E} $. 
\end{proof}

\bigskip
Combining Corollaries 
\ref{quantile-order-ES-c} and 
\ref{quantile-finite-c} yields a condition under which   $ \LL_Q = \LL_{Q^E} $   and   $ \varrho_Q $   is finite. 
Corollary 
\ref{quantile-order-1-c} also yields some further results on the comparison of quantile risk measures and their domains: 

\bigskip
\begin{corollary}{}
\label{quantile-order-2-c}
\begin{onelist}
\item   If   $ D_1 \leq D_2, $ 
then   $ \LL_{Q_1} \subseteq \LL_{Q_2} $   and   $ \varrho_{Q_2}[X] \leq \varrho_{Q_1}[X] $   holds for every   $ X \in \LL_{Q_1} $. 
\item   If   $ D \leq D^E, $ 
then   $ \LL_Q \subseteq \LL_{Q^E} $   and   $ E[X] \leq \varrho_Q[X] $   holds for every   $ X \in \LL_Q $. 
\item   If   $ D $   is convex, 
then   $ \LL_Q \subseteq \LL_{Q^E} $   and   $ E[X] \leq \varrho_Q[X] $   holds for every   $ X \in \LL_Q $. 
\item   Consider   $ \alpha,\beta\in[0,1) $. 
If   $ \alpha\leq\beta, $ 
then   $ \varrho_{Q^{\ES_\alpha}}[X] \leq \varrho_{Q^{\ES_\beta}}[X] $   holds for every   $ X \in \LL_{Q^E} $. 
\item   The identity   $ E[X] = \inf_{\alpha\in(0,1)} \varrho_{Q^{\ES_\alpha}}[X] $   holds for every   $ X \in \LL_{Q^E} $. 
\end{onelist}
\end{corollary}

\bigskip
\begin{proof}
Assertion (1) is immediate from Theorem 
\ref{quantile-char-t} and Corollary 
\ref{quantile-order-1-c} and yields assertions (2), (3) and (4). 
Assertion (5) follows from the dominated convergence theorem. 
\end{proof}

%\bigskip
%Because of the particular role of the set   $ \LL_{Q^E} $, 
%we note the following result: 
%
%\bigskip
%\begin{lemma}{}
%\label{quantile-E-cone-l}
%$ \LL_{Q^E} $   is a convex cone and   $ \varrho_{Q^E} $   is subadditive. 
%\end{lemma}
%
%\bigskip
%\begin{proof}
%Consider   $ X,Y\in\LL^0 $   and   $ a\in\R_+ $. 
%Then we have   $ E[(aX)^+] = a\,E[X^+] $   and   $ E[(X\!+\!Y)^+] \leq  E[X^+] + E[Y^+] $, 
%which implies that   $ \LL_{Q^E} $   is a convex cone. 
%Moreover, 
%we have   $ \varrho_{Q^E}[X\!+\!Y] = E[X\!+\!Y] = E[X] + E[Y] = \varrho_{Q^E}[X] + \varrho_{Q^E}[Y] $, 
%which implies that   $ \varrho_{Q^E} $   is subadditive. 
%\end{proof}
%
%\bigskip
%On the other hand, 
%the following example shows that   $ \LL_{Q^E} $   fails to be a vector space: 
%
%\bigskip
%\begin{example}{}
%\label{quantile-E-e}
%Consider a random variable   $ X $   satisfying 
%\begin{eqnarray*}
%        P_X                                                                   
%&  = &  \int \frac{2}{\pi}\,\frac{1}{1+x^2}\,\chi_{(-\infty,0)}(x)\,d\leb(x)  
%\end{eqnarray*}
%Then we have   $ E[X^+] = 0 $    and   $ E[(-X)^+] = E[X^-] =\infty $, 
%and hence   $ X\in\LL_{Q^E} $   and   $ -X\notin\LL_{Q^E} $. 
%\end{example}

%%%%%%%%%%%%%%%%%%%%%%%%%%%%%%%%%%%%%%%%%%%%%%%%%%%%%%%%%%%%%%%%%%%%%%%%%%%%%%%%%%%%%%%%%%%%%%%%%%%

\section{Spectral Risk Measures}
\label{spectral}

A map   $ s : (0,1)\to\R_+ $   is said to be a 
\emph{spectral function} if it is increasing and satisfies   $ \int_{(0,1)} s(u)\,d\leb(u) = 1 $. 

\bigskip
The quantile risk measure   $ \varrho_Q $   is said to be a 
\emph{spectral risk measure} if there exists a spectral function   $ s $   such that 
\begin{eqnarray*}
        Q                    
&  = &  \int s(u)\,d\leb(u)  
\end{eqnarray*}
Thus, 
if   $ \varrho_Q $   is a spectral risk measure with spectral function   $ s $, 
then the identities 
\begin{eqnarray*}
        \LL_Q                                                                                            
&  = &  \biggl\{ X\in\LL^0 \biggm| \int_{(0,1)} (F^\leftarrow_X(u))^+\,s(u)\,d\leb(u) < \infty \biggr\}  
\end{eqnarray*}
and 
\begin{eqnarray*}
        \varrho_Q[X]                                    
&  = &  \int_{(0,1)} F^\leftarrow_X(u)\,s(u)\,d\leb(u)  
\end{eqnarray*}
hold for every   $ X \in \LL_Q $. 
Note that the spectral function of a spectral risk measure is unique almost everywhere, 
by the Radon--Nikodym theorem. 

\bigskip
\begin{examples}{}
\label{spectral-e}
\begin{onelist}
\item   {\bf Expectation:} 
Since   $ D^E(u) = u $, 
we have 
\begin{eqnarray*}
        Q^E   
&  = &  \leb  
\end{eqnarray*}
and the function   $ s^E : (0,1)\to\R_+ $   given by 
\begin{eqnarray*}
        s^E(u)  
& := &  1       
\end{eqnarray*}
is a spectral function. 
Therefore, 
$ \varrho_{Q^E} $   is a spectral risk measure. 
\item   {\bf Value at Risk:}
For every   $ \alpha\in(0,1) $, 
$ Q^{\VaR_\alpha} $   is the Dirac measure at   $ \alpha $    and hence does not have a density with respect to   $ \leb $. 
Therefore, 
$ \varrho_{Q^{\VaR_\alpha}} $   is not a spectral risk measure. 
\item   {\bf Expected Shortfall:}
For every   $ \alpha\in[0,1) $,
we have 
\begin{eqnarray*}
        Q^{\ES_\alpha}                                           
&  = &  \int \frac{1}{1-\alpha}\,\chi_{(\alpha,1)}(u)\,d\leb(u)  
\end{eqnarray*}
and the function   $ s^{\ES_\alpha} : (0,1)\to\R_+ $   given by 
\begin{eqnarray*}
        s^{\ES_\alpha}(u)                         
& := &  \frac{1}{1-\alpha}\,\chi_{(\alpha,1)}(u)  
\end{eqnarray*}
is a spectral function. 
Therefore, 
$ \varrho_{Q^{\ES_\alpha}} $   is a spectral risk measure. 
\item   {\bf Expected Shortfall of Higher Order:} 
For every   $ n\in\N $   and   $ \alpha\in[0,1) $, 
we have 
\begin{eqnarray*}
        Q^{\ES_{n,\alpha}}                                                                                       
&  = &  \int \frac{n}{1-\alpha} \biggl( \frac{u-\alpha}{1-\alpha} \biggr)^{n-1}\,\chi_{(\alpha,1)}(u)\,d\leb(u)  
\end{eqnarray*}
and the function   $ s^{\ES_{n,\alpha}} : (0,1)\to\R_+ $   given by 
\begin{eqnarray*}
        s^{\ES_{n,\alpha}}(u)                                                                     
& := &  \frac{n}{1-\alpha} \biggl( \frac{u-\alpha}{1-\alpha} \biggr)^{n-1}\,\chi_{(\alpha,1)}(u)  
\end{eqnarray*}
is a spectral function. 
Therefore, 
$ \varrho_{Q^{\ES_{n,\alpha}}} $   is a spectral risk measure. 
\end{onelist}
\end{examples}

\bigskip
Our aim is to characterize the spectral risk measures within the class of all quantile risk measures. 
The following result is inspired by 
Gzyl and Mayoral [2008] who considered distortion risk measures on the positive cone of   $ \LL^2 $: 

\bigskip
\begin{theorem}{}
\label{spectral-char-t}
The following are equivalent:
\begin{abclist}
\item   $ D $   is convex. 
\item   There exists a spectral function   $ s $   such that   $ Q = \int s(u)\,d\leb(u) $. 
\item   $ \varrho_Q $   is a spectral risk measure. 
\end{abclist}
In this case, 
every spectral function   $ s $   representing   $ Q $   satisfies   $ s = D^\one $   almost everywhere $($with respect to   $ \lambda) $. 
\end{theorem}

\bigskip
\begin{proof}
Since   $ \lim_{u \to 0} D(u) = 0 = D(0) $   and   $ \lim_{u \to 1} D(u) = 1 = D(1) $, 
$ D $   is convex if and only if   $ D $   is convex on   $ (0,1) $. 
\\
Assume first that (a) holds. 
The following arguments are taken from Aliprantis and Burkinshaw [1990; Chapter 29]. 
Since   $ D $   is increasing, 
$ D $   is differentiable almost everywhere, 
and 
since   $ D $   is convex, 
its derivative   $ D^\one $   is increasing. 
Consider now an arbitrary interval   $ [u,v]\subseteq(0,1) $. 
Since   $ D $   is convex, 
the restriction of   $ D $   to   $ [u,v] $   is Lipschitz continuous, 
hence absolutely continuous, 
and thus continuous and of bounded variation.  
Therefore, 
the restriction of   $ Q $   to the $\sigma$--field of all Borel sets in   $ [u,v] $   
is absolutely continuous with respect to the restriction of   $ \lambda $   and its Radon--Nikodym derivative agrees with   $ D^\one $. 
Since   $ [u,v]\subseteq(0,1) $   was arbitrary, 
it follows that   $ Q $   is absolutely continuous with respect to   $ \lambda $, 
and 
since the Radon--Nikodym derivative   $ s : (0,1)\to\R_+ $   of   $ Q $   with respect to   $ \leb $   is unique almost everywhere, 
it follows that   $ s = D^\one $   almost everywhere. 
This yields the existence of an increasing function   $ s : (0,1)\to\R_+ $   satisfying   $ Q = \int s(u)\,d\leb(u) $. 
%\begin{eqnarray*}
%        Q                    
%&  = &  \int s(u)\,d\leb(u)  
%\end{eqnarray*}
Therefore, 
(a) implies (b). 
\\
Assume now that (b) holds. 
Since   $ s $   is increasing, 
we have, 
for any   $ u,v,w\in(0,1) $   such that   $ u < v < w $, 
\begin{eqnarray*}
        \frac{D(v)-D(u)}{v-u}                      
&  = &  \frac{1}{v-u} \int_{(u,v]} s(t)\,d\leb(t)  \\[1ex]
&\leq&  s(v)                                       \\[1ex]
&\leq&  \frac{1}{w-v} \int_{(v,w]} s(t)\,d\leb(t)  \\*
&  = &  \frac{D(w)-D(v)}{w-v}                      
\end{eqnarray*}
which implies that   $ D $   is convex. 
Therefore, 
implies (a).  
\end{proof}

\bigskip
\begin{theorem}{}
\label{spectral-mixture-t}
If   $ D $   is convex, 
then there exists a measure   $ \nu : \BB([0,1))\to[0,\infty] $   such that 
$$  \varrho_Q[X]  =  \int_{[0,1)} (1\!-\!\alpha)\,\varrho_{Q^{\ES_\alpha}}[X]\,d\nu(\alpha)  $$
holds for every   $ X\in\LL_Q $. 
\end{theorem}
\pagebreak

\bigskip
\begin{proof}
Without loss of generality, 
we may and do assume that   $ s $   is continuous from the right. 
Define   $ s(0) := \inf_{u\in(0,1)} s(u) $. 
Then there exists a unique $\sigma$--finite measure   $ \nu : \BB([0,1))\to[0,\infty] $   satisfying   $ \nu[[0,u]] = s(u) $   for all   $ u\in(0,1) $. 
Since the map   $ (0,1)\times[0,1)\to\R : (u,\alpha) \to F_X^\leftarrow(u)\,\chi_{[0,u]}(\alpha) $   is measurable 
and its positive part is integrable with respect to the product measure   $ \nu\otimes\leb $, 
Fubini's theorem yields 
\begin{eqnarray*}
        \varrho_Q[X]                                                                               
&  = &  \int_{(0,1)} F^\leftarrow_X(u)\,s(u)                                           \,d\leb(u)  \\
&  = &  \int_{(0,1)} F^\leftarrow_X(u) \int_{[0,1)}  \chi_{[0,u]}(\alpha)\,d\nu(\alpha)\,d\leb(u)  \\
&  = &  \int_{[0,1)} \int_{(0,1)} F^\leftarrow_X(u)\,\chi_{(\alpha,1)}(u)\,d\leb(u)\,d\nu(\alpha)  \\*
&  = &  \int_{[0,1)} (1\!-\!\alpha)\,\varrho_{Q^{\ES_\alpha}}[X]                   \,d\nu(\alpha)  
\end{eqnarray*}
This proves the assertion. 
\end{proof}

%%%%%%%%%%%%%%%%%%%%%%%%%%%%%%%%%%%%%%%%%%%%%%%%%%%%%%%%%%%%%%%%%%%%%%%%%%%%%%%%%%%%%%%%%%%%%%%%%%%

\section{Subadditivity of Spectral Risk Measures}
\label{subadditive}

In the present section we show that a quantile risk measure is subadditive if and only if its distortion function is convex. 
To prove that convexity of the distortion function is sufficient for subadditivity of the quantile risk measure, 
we use Theorem 
\ref{spectral-mixture-t}. 
Since the expectation is additive and hence subadditive, 
it remains to show that the expected shortfall at any level is subadditive. 

\bigskip
To establish subadditivity of the expected shortfall we need the following lemma 
from which we deduce another representation of the values of the expected shortfall: 

\bigskip
\begin{lemma}{}
\label{subadditive-SL-l}
The identity 
$$  \int_{(0,1)} \Bigl( F_X^\leftarrow(u) - c \Bigr)^+\,d\leb(u)  =  E[(X-c)^+]  $$
holds for every   $ X\in\LL^0 $   and   $ c\in\R $. 
\end{lemma}

\bigskip
\begin{proof}
Lemma 
\ref{preliminaries-l} together with Lemma 
\ref{quantile-char-l} yields 
\begin{eqnarray*}
        \int_{(0,1)} \Bigl( F_X^\leftarrow(u) - c \Bigr)^+\,d\leb(u)  
&  = &  \int_{(0,1)}   (F_{X-c}^\leftarrow(u))^+          \,d\leb(u)  \\
&  = &  \int_\R z^+    \,dF_{X-c}(z)                                  \\
&  = &  \int_\R (x-c)^+\,dF_X(x)                                      \\*[1ex]
&  = &  E[(X-c)^+]                                                    
\end{eqnarray*}
as was to be shown. 
\end{proof}

\bigskip
\begin{lemma}{}
\label{subadditive-ES-l}
For every   $ \alpha\in(0,1), $ 
the identity 
$$  \varrho_{Q^{\ES_\alpha}}[X]                                                                          
 =  F_X^\leftarrow(\alpha) + \frac{1}{1-\alpha}\,E\Bigl[ \Bigl(X-F_X^\leftarrow(\alpha) \Bigr)^+ \Bigr]  
 =  \inf_{c\in\R} \biggl( c + \frac{1}{1-\alpha}\,E[(X\!-\!c)^+] \biggr)                                 
$$
holds for every   $ X\in\LL_{Q^{\ES_\alpha}} $. 
\end{lemma}

\bigskip
\begin{proof}
Because of Lemma 
\ref{subadditive-SL-l}, 
we have 
\begin{eqnarray*}
        \varrho_{Q^{\ES_\alpha}}[X]                                                        - F_X^\leftarrow(\alpha)  
&  = &  \int_{(0,1)} F_X^\leftarrow(u)\,\frac{1}{1-\alpha}\,\chi_{(\alpha,1)}(u)\,d\leb(u) - F_X^\leftarrow(\alpha)  \\
&  = &  \frac{1}{1-\alpha} \int_{(0,1)} \Bigl( F_X^\leftarrow(u) - F_X^\leftarrow(\alpha) \Bigr)^+\,d\leb(u)         \\*
&  = &  \frac{1}{1-\alpha}\,E\Bigl[ \Bigl(X-F_X^\leftarrow(\alpha) \Bigr)^+ \Bigr]                                   
\end{eqnarray*}
which yields the first identity for   $ \varrho_{Q^{\ES_\alpha}}[X] $. 
Now the second identity for   $ \varrho_{Q^{\ES_\alpha}}[X] $   is easily verified 
by a distinction of the cases   $ c<F_X^\leftarrow(\alpha) $   and   $ c>F_X^\leftarrow(\alpha) $. 
\end{proof}

\bigskip
Lemmas 
\ref{subadditive-SL-l} and 
\ref{subadditive-ES-l} are well--known and are frequently used to establish the subadditivity of expected shortfall on   $ \LL^\infty $; 
see e.\,g.\ 
Embrechts and Wang [2015]. 
The following result also relies on these lemmas, 
but it is more general and more precise: 

\bigskip
\begin{lemma}{}
\label{subadditive-ES-cone-l}
For every   $ \alpha\in[0,1) $, 
$ \LL_{Q^{\ES_\alpha}} $   is a convex cone and $ \varrho_{Q^{\ES_\alpha}} $   is subadditive. 
\end{lemma}

\bigskip
\begin{proof}
%By Example 
%\ref{quantile-e}(3), 
Since   $ \LL_{Q^{\ES_\alpha}}=\LL_{Q^E} $, 
we see that 
%Lemma 
%\ref{quantile-E-cone-l} implies that   
$ \LL_{Q^{\ES_\alpha}} $   is a convex cone. 
%Moreover, 
Also, 
since   $ Q^{\ES_0}=Q^E $, 
%Lemma 
%\ref{quantile-E-cone-l} implies that   
we see that 
$ \varrho_{Q^{\ES_0}} $   is subadditive. 
Consider now   $ \alpha\in(0,1) $   and   $ X,Y\in\LL_{Q^{\ES_\alpha}} $. 
Then we have   $ X+Y\in\LL_{Q^{\ES_\alpha}} $   and, 
for any   $ x,y\in\R $, 
Lemma 
\ref{subadditive-ES-l} yields 
\begin{eqnarray*}
        \varrho_{Q^{\ES_\alpha}}[X+Y]                                             
&\leq&  (x+y) + \frac{1}{1-\alpha}\,E\Bigl[ \Bigl( (X+Y) - (x+y) \Bigr)^+ \Bigr]  \\*
&  = &   x+y  + \frac{1}{1-\alpha}\,E\Bigl[ \Bigl( (X-x) + (Y-y) \Bigr)^+ \Bigr]  
\end{eqnarray*}
and hence 
\begin{eqnarray*}
        \varrho_{Q^{\ES_\alpha}}[X+Y]                                                                                    
&\leq&  \biggl( x + \frac{1}{1-\alpha}\,E[(X\!-\!x)^+] \biggr) + \biggl( y + \frac{1}{1-\alpha}\,E[(Y\!-\!y)^+] \biggr)  
\end{eqnarray*}
Now minimization over   $ x,y\in\R $   and using Lemma 
\ref{subadditive-ES-l} again yields 
\begin{eqnarray*}
        \varrho_{Q^{\ES_\alpha}}[X+Y]                              
&\leq&  \varrho_{Q^{\ES_\alpha}}[X] + \varrho_{Q^{\ES_\alpha}}[Y]  
\end{eqnarray*}
Therefore, 
$ \varrho_{Q^{\ES_\alpha}} $   is subadditive for every   $ \alpha\in(0,1) $. 
\end{proof}

\bigskip
The previous result provides the key for proving the main implication of the following theorem; 
see also 
Wang and Dhaene [1998] who considered distortion risk measures on the positive cone of   $ \LL^1 $   
and used a proof based on comonotonicity. 
\pagebreak

\bigskip
\begin{theorem}{}
\label{subadditive-char-t}
The following are equivalent: 
\begin{abclist}
\item   $ D $   is convex. 
\item   $ \varrho_Q $   is subadditive. 
\item   $ \LL_Q $   is a convex cone and   $ \varrho_Q $   is subadditive. 
\end{abclist}%
\end{theorem}

\bigskip
\begin{proof}
Assume first that (a) holds 
and  
%Then Theorem 
%\ref{spectral-mixture-t} together with Lemma 
%\ref{subadditive-ES-cone-l} implies that   $ \varrho_Q $   is subadditive. 
consider a spectral function   $ s $   representing   $ Q $   and the measure   $ \nu $   constructed in the proof of Theorem 
\ref{spectral-mixture-t}. 
Consider   $ X,Y\in\LL_Q $   and   $ a\in\R_+ $. 
Then we have   $ aX\in\LL_Q $. 
Moreover,  
since   $ D $   is convex, 
Corollary 
\ref{quantile-order-2-c} yields   $ X,Y\in\LL_{Q^E} $.   
For every   $ \alpha\in[0,1) $, 
this yields   $ X,Y\in\LL_{Q^{\ES_\alpha}} $, 
hence   $ X+Y\in\LL_{Q^{\ES_\alpha}} $, 
by Lemma 
\ref{subadditive-ES-cone-l}, 
and thus   $ X^+,Y^+,(X+Y)^+\in\LL_{Q^{\ES_\alpha}} $. 
Proceeding as in the proof of Theorem 
\ref{spectral-mixture-t} and using Lemma 
\ref{subadditive-ES-cone-l} again, 
we obtain 
\begin{eqnarray*} 
        \int_{(0,1)} F^\leftarrow_{(X+Y)^+}(u)\,s(u)                                                            \,d\leb(u)      
&  = &  \int_{[0,1)} (1\!-\!\alpha)\,\varrho_{Q^{\ES_\alpha}}[(X+Y)^+]                                          \,d\nu(\alpha)  \\
&\leq&  \int_{[0,1)} (1\!-\!\alpha)\,\Bigl( \varrho_{Q^{\ES_\alpha}}[X^+] + \varrho_{Q^{\ES_\alpha}}[Y^+] \Bigr)\,d\nu(\alpha)  \\
&  = &  \int_{[0,1)} (1\!-\!\alpha)\,\varrho_{Q^{\ES_\alpha}}[X^+]                                              \,d\nu(\alpha)  
      + \int_{[0,1)} (1\!-\!\alpha)\,\varrho_{Q^{\ES_\alpha}}[Y^+]                                              \,d\nu(\alpha)  \\[.5ex]
&  = &  \varrho_Q[X^+] + \varrho_Q[Y^+]                                                                                         \\*[1ex]
&  < &  \infty                                                                                                                  
\end{eqnarray*}
This yields   $ (X+Y)^+\in\LL_Q $   and hence   $ X+Y\in\LL_Q $. 
Thus, 
$ \LL_Q $   is a convex cone, 
and Theorem 
\ref{spectral-mixture-t} together with Lemma 
\ref{subadditive-ES-cone-l} implies that   $ \varrho_Q $   is subadditive. 
Therefore, 
(a) implies (c). 
Obviously, 
(c) implies (b), 
and it follows from Example 
\ref{subadditive-e} below that 
(b) implies (a).  
\end{proof}

\bigskip
For the discussion of the subsequent Example 
\ref{subadditive-e} we need the following lemma: 

\bigskip
\begin{lemma}{}
\label{subadditive-convexity-l}
The following are equivalent: 
\begin{abclist}
\item   $ D $   is convex. 
\item   The inequality   
$$  D(u)  \leq  \frac{1}{2} \Bigl( D(u-\varepsilon) + D(u+\varepsilon) \Bigr)  $$
holds for all   $ u\in(0,1) $   and   $ \varepsilon\in(0,\min\{u,1\!-\!u\}) $. 
\end{abclist}
\end{lemma}

\bigskip
\begin{proof}
Assume that (b) holds. 
Then the inequality 
\begin{eqnarray*}
        D\biggl( \frac{u+v}{2} \biggr)         
&\leq&  \frac{1}{2} \Bigl( D(u) + D(v) \Bigr)  
\end{eqnarray*}
holds for all   $ u,v\in(0,1) $. 
Furthermore, 
since   $ D $   is monotone, 
the left and right limits   $ D(u+) $   and   $ D(u-) $   exist for every   $ u\in(0,1) $. 
From   $ 2\,D(u) \leq D(u-\varepsilon) + D(u+\varepsilon) $   we obtain 
\begin{eqnarray*}
        2\,D(u)        
&\leq&  D(u-) + D(u+)  
\end{eqnarray*}
and from   $ 2\,D(u+\varepsilon) \leq D(u) + D(u+2\,\varepsilon) $   we obtain 
\begin{eqnarray*}
        D(u+)  
&\leq&  D(u)   
\end{eqnarray*}
Combining these two inequalities yields   $ D(u) \leq D(u-) \leq D(u+) \leq D(u)  $,   
and hence   $ D(u-) = D(u) = D(u+) $. 
Therefore, 
$ D $   is continuous, 
and it now follows from the first inequality of this proof that   $ D $   is convex. 
Therefore, 
(b) implies (a). 
The converse implication is obvious. 
\end{proof}

\bigskip
The bivariate distribution discussed in the following example was proposed by 
Wirch and Hardy [2002]. 

\bigskip
\begin{example}{}
\label{subadditive-e}
Assume that   $ D $   is not convex. 
Then Lemma 
\ref{subadditive-convexity-l} yields the existence of some   $ u\in(0,1) $   and   $ \varepsilon\in(0,\min\{u,1\!-\!u\}) $   such that 
\begin{eqnarray*}
        2\,D(u)                              
&  > &  D(u-\varepsilon) + D(u+\varepsilon)  
\end{eqnarray*}
Consider random variables   $ X,Y\in\LL^\infty $   whose joint distribution is given by the following table with   $ a\in(0,\infty) $: 
$$\begin{array}{|@{\;}c@{\;}|@{\;}c@{\quad}c@{\quad}c@{\;}|@{\;}c@{\;}@{\;}c@{\;}|}
\hline\rule{0ex}{2.5ex}%
                   & \multicolumn{3}{c|@{\;}}{y}  & & \\
 x	               & -(a+\varepsilon) & -(a+\varepsilon/2) & 0     & P[\{X=x\}] & P[\{X \leq x\}] \\
\hline\rule{0ex}{2.5ex}%
 -(a+\varepsilon)  & u-\varepsilon    & 0                  &     \varepsilon    &   u        & u  \\
 0                 & 0                & \varepsilon        & 1-u-\varepsilon    & 1-u        & 1  \\
\hline\rule{0ex}{2.5ex}%
P[\{Y=y\}]         & u-\varepsilon    & \varepsilon        & 1-u   &            &                 \\
\rule{0ex}{2.5ex}%
P[\{Y \leq y\}]    & u-\varepsilon    & u                  & 1     &            &                 \\
\hline
\end{array}$$
Then the distribution of the sum   $ X+Y $   is given by the table 
$$\begin{array}{|@{\;}c@{\;}|@{\;}c@{\quad}c@{\quad}c@{\quad}@{\;}c@{\;}|}
\hline\rule{0ex}{2.5ex}%
 z	               & -2(a+\varepsilon) & -(a+\varepsilon) & -(a+\varepsilon/2) & 0                 \\
\hline\rule{0ex}{2.5ex}%
P[\{X+Y=z\}]       &    u-\varepsilon    & \varepsilon    &   \varepsilon      & 1-u-\varepsilon   \\
\rule{0ex}{2.5ex}%
P[\{X+Y \leq z\}]  &    u-\varepsilon    & u              & u+\varepsilon      & 1                 \\
\hline
\end{array}$$
Because of Theorem 
\ref{quantile-char-t} this yields 
\begin{eqnarray*}
        \varrho_Q[X]                                                                                                
&  = &  -\,(a+\varepsilon)\,D(u)                                                                                    \\
        \varrho_Q[Y]                                                                                                
&  = &  -\,(\varepsilon/2)\,D(u\!-\!\varepsilon) - (a+\varepsilon/2)\,D(u)                                          \\*
        \varrho_Q[X\!+\!Y]                                                                                          
&  = &  -\,(a+\varepsilon)\,D(u\!-\!\varepsilon) - (\varepsilon/2)\,D(u) - (a+\varepsilon/2)\,D(u\!+\!\varepsilon)  
\end{eqnarray*}
and hence 
\begin{eqnarray*}
        \varrho_Q[X\!+\!Y]                                                                                                    
&  = &  \varrho_Q[X] + \varrho_Q[Y] + (a+\varepsilon/2)\,\Bigl( 2\,D(u) - D(u\!-\!\varepsilon) - D(u\!+\!\varepsilon) \Bigr)  \\*
&  > &  \varrho_Q[X] + \varrho_Q[Y]                                                                                           
\end{eqnarray*}
Therefore, 
$ \varrho_Q $   fails to be subadditive. 
\end{example}

%%%%%%%%%%%%%%%%%%%%%%%%%%%%%%%%%%%%%%%%%%%%%%%%%%%%%%%%%%%%%%%%%%%%%%%%%%%%%%%%%%%%%%%%%%%%%%%%%%%

\section{On the Domain of a Quantile Risk Measure}
\label{domain}

In this final section we compare the domain 
\begin{eqnarray*}
        \LL_Q                                                                                   
&  = &  \biggl\{ X\in\LL^0 \biggm| \int_{(0,1)} (F^\leftarrow_X(u))^+\,dQ(u) < \infty \biggr\}  
\end{eqnarray*}
of the quantile risk measure   $ \varrho_Q $   with two other classes of random variables. 
Define 
\begin{eqnarray*}
        \LL_Q^\Acerbi                                                                           
& := &  \biggl\{ X\in\LL^0 \biggm| \int_{(0,1)} |F^\leftarrow_X(u)|  \,dQ(u) < \infty \biggr\}  
\end{eqnarray*}
and 
\begin{eqnarray*}
        \LL_Q^\Pichler                                                                          
& := &  \biggl\{ X\in\LL^0 \biggm| \int_{(0,1)} F^\leftarrow_{|X|}(u)\,dQ(u) < \infty \biggr\}  
\end{eqnarray*}
In the case where   $ Q $   is represented by a spectral function, 
these classes were introduced by 
Acerbi [2002] and 
Pichler [2013], 
respectively. 
We have   $ \LL_Q^\Acerbi\subseteq\LL_Q $, 
and Corollary 
\ref{quantile-finite-c} provides a sufficient condition for   $ \LL_Q^\Acerbi=\LL_Q $. 
Moreover, 
since   $ X^+ \leq |X| $, 
we also have   $ \LL_Q^\Pichler\subseteq\LL_Q $. 
Below we shall show that   $ \LL_Q^\Pichler\subseteq\LL_Q^\Acerbi $   whenever   $ D $   is convex. 
To this end, 
we need the following lemma: 

\bigskip
\begin{lemma}{}
\label{domain-A-l}
Assume that   $ D $   is convex and consider   $ X\in\LL^0 $. 
If   $ X^+\in\LL_Q^\Acerbi $   and   $ X^-\in\LL_Q^\Acerbi, $ 
then   $ X\in\LL_Q^\Acerbi $. 
\end{lemma}

\bigskip
\begin{proof}
From   $ (F_X^\leftarrow)^+ = F_{X^+}^\leftarrow $   and   $ X^+\in\LL_Q^\Acerbi $   we obtain 
\begin{eqnarray*}
        \int_{(0,1)} (F_X^\leftarrow(u))^+\,dQ(u)  
&  < &  \infty                                     
\end{eqnarray*}
To prove that the integral   $ \int_{(0,1)} (F_X^\leftarrow(u))^-\,dQ(u) $   is finite as well, 
we need the upper quantile function   $ F_X^\rightarrow: (0,1) \to \R $   given by 
\begin{eqnarray*}
        F_X^\rightarrow(u)                                
& := &  \sup\Bigl\{ x\in\R \Bigm| F_X (x) \leq u \Bigr\}  
\end{eqnarray*}
The lower and upper quantile functions satisfy   $ F_X^\leftarrow \leq F_X^\rightarrow $, 
and we have 
\begin{eqnarray*}
        (F_X^\leftarrow(u))^-                          
&  = &  -\,F_X^\leftarrow (u)\,\chi_{(0,F_{X}(0)]}(u)  
\end{eqnarray*}
and 
\begin{eqnarray*}
        F_{X^-}^\leftarrow(1\!-\!u)                    
&  = &  -\,F_X^\rightarrow(u)\,\chi_{(0,F_{X}(0))}(u)  
\end{eqnarray*}
almost everywhere with respect to   $ \leb $. 
Since   $ D $   is convex and hence continuous, 
$ Q $   is absolutely continuous with respect to   $ \leb $. 
This yields 
\begin{eqnarray*}
        0                                                                                          
&\leq&  \int_{(0,1)}        \Bigl(  F_X^\rightarrow(u) - F_X^\leftarrow (u) \Bigr)        \,dQ(u)  \\
&  = &  \int_{(0,1)} \int_\R \chi_{[F_X^\leftarrow (u) , F_X^\rightarrow(u))}(x)\,d\leb(x)\,dQ(u)  \\
&\leq&  \int_\R \int_{(0,1)} \chi_{\{F_{X}(x)\}}   (u)                          \,dQ(u)\,d\leb(x)  \\*[1ex]
&  = &  0                                                                                          
\end{eqnarray*}
and hence   $ F_X^\rightarrow = F_X^\leftarrow $   almost everywhere with respect to   $ Q $. 
Consider now a spectral function   $ s $   representing   $ Q $. 
Since   $ s $   is positive and increasing, 
we obtain 
\begin{eqnarray*}
        \int_{(0,1)} (F_X^\leftarrow(u))^-                       \,dQ(u)                 
&  = &  \int_{(0,1)} (- F_X^\leftarrow (u))\,\chi_{(0,F_X(0)]}(u)\,dQ(u)                 \\
&  = &  \int_{(0,1)} (- F_X^\rightarrow(u))\,\chi_{(0,F_X(0))}(u)\,dQ(u)                 \\*
&  = &  \int_{(0,1)} (- F_X^\rightarrow(u))\,\chi_{(0,F_X(0))}(u)\,s(u)      \,d\leb(u)  \\*
&  = &  \int_{(0,1)}   F_{X^-}^\leftarrow(1\!-\!u)               \,s(u)      \,d\leb(u)  \\
&  = &  \int_{(0,1)}   F_{X^-}^\leftarrow(u)                     \,s(1\!-\!u)\,d\leb(u)  \\
&\leq&  \int_{(0,1/2)} F_{X^-}^\leftarrow(1/2)                   \,s(1\!-\!u)\,d\leb(u)  
      + \int_{(1/2,1)} F_{X^-}^\leftarrow(u)                     \,s(u)      \,d\leb(u)  \\*
&\leq&                 F_{X^-}^\leftarrow(1/2)                                           
      + \int_{(0,1)}   F_{X^-}^\leftarrow(u)                     \,dQ(u)                 
\end{eqnarray*}
Since   $ X^-\in\LL_Q^\Acerbi $, 
the last expression is finite, 
and this yields 
\begin{eqnarray*}
        \int_{(0,1)} (F_X^\leftarrow(u))^-\,dQ(u)  
&  < &  \infty                                     
\end{eqnarray*}
Therefore, 
we have   $ X\in\LL_Q^\Acerbi $. 
\end{proof}

\bigskip
\begin{theorem}{}
\label{domain-A-t}
If   $ D $   is convex, 
then   $ \LL_Q^\Pichler\subseteq\LL_Q^\Acerbi $
\end{theorem}

\bigskip
\begin{proof}
Consider   $ X\in\LL_Q^\Pichler $. 
Then we have   $ |X|\in\LL_Q^\Pichler $, 
hence   $ X^+,X^-\in\LL_Q^\Pichler $, 
and thus   $ X^+,X^-\in\LL_Q^\Acerbi $. 
Now Lemma 
\ref{domain-A-l} yields   $ X\in\LL_Q^\Acerbi $. 
\end{proof}

\bigskip
The following examples provide some further insight into the relationships between these three classes of random variables: 

\bigskip
\begin{examples}{}
\label{comain-e}
\begin{onelist}
\item   If   $ D = D^{\VaR_\alpha} $, 
then   $ \LL_Q^\Pichler = \LL_Q^\Acerbi = \LL_Q = \LL^0 $. 
\item   If   $ D = D^E $, 
then   $ \LL_Q^\Pichler = \LL_Q^\Acerbi = \LL^1 \neq\LL_Q $. 
\item   If   $ D = D^{\ES_\alpha} $   for some   $ \alpha\in(0,1) $, 
then   $ \LL_Q^\Pichler \neq \LL^Q = \LL_Q^\Acerbi $. 
\item   Assume that there exists some   $ \delta\in(0,1) $   such that   $ D $   satisfies 
\begin{eqnarray*}
        D(u)                                            
&  = &  u\,\chi_{[0,\delta)}(u) + \chi_{[\delta,1]}(u)  
\end{eqnarray*}
(and hence fails to be convex). 
Then every   $ X\in\LL^0 $   satisfies   
$$  \int_{(0,1)} F_{|X|}^\leftarrow(u)\,dQ(u)  <  \infty  \qquad\text{and}\qquad  \int_{(0,1)} (F_X^\leftarrow(u))^+\,dQ(u)  <  \infty  $$
This yields   $  \LL_Q^\Pichler  =  \LL^0  =  \LL_Q  $, 
as well as 
\begin{eqnarray*}
        \LL_Q^\Acerbi                                                                                   
&  = &  \biggl\{ X\in\LL^0 \biggm| \int_{(0,1)}      (F_X^\leftarrow(u))^-\,dQ   (u) < \infty \biggr\}  \\
&  = &  \biggl\{ X\in\LL^0 \biggm| \int_{(0,\delta)} (F_X^\leftarrow(u))^-\,d\leb(u) < \infty \biggr\}  \\
&  = &  \biggl\{ X\in\LL^0 \biggm| \int_{(0,1)}      (F_X^\leftarrow(u))^-\,d\leb(u) < \infty \biggr\}  \\*[.5ex]
&  = &  \Bigl \{ X\in\LL^0 \Bigm | E[X^-]                                            < \infty \Bigr \}  
\end{eqnarray*}
such that   $ \LL_Q^\Pichler \neq \LL_Q^\Acerbi $   and   $ \LL_Q^\Acerbi \neq \LL_Q $. 
\item   Assume that   $ D $   satisfies 
\begin{eqnarray*}
        D(u)                                                             
&  = &  \frac{1}{2}\,\sqrt{u}\,\chi_{[0,1/4)}(u) + u\,\chi_{[1/4,1]}(u)  
\end{eqnarray*}
%(and hence fails to be convex). 
Then Corollary 
\ref{quantile-order-ES-c} yields   $ \LL_Q = \LL_{Q^E} $. 
Moreover, 
straightforward calculation yields 
\begin{eqnarray*}
        \int_{(0,1)} F_{|X|}^\leftarrow(u)\,dQ(u) 
&\leq&  \leb[(0,F_{|X|}^\leftarrow(1/4))]
%\int_{(0,F_{|X|}^\leftarrow(1/4))}      \Bigl( 1 - (D \circ F_{|X|})(0) \Bigr)\,d\leb(x)  \\*
   + \int_{[F_{|X|}^\leftarrow(1/4),\infty)} \Bigl( 1 - (D \circ F_{|X|})(x) \Bigr)\,d\leb(x) 
\end{eqnarray*}
and 
\begin{eqnarray*}
        \int_{(0,1)} F_{|X|}^\leftarrow(u)\,d\leb(u) 
&\leq&  \leb[(0,F_{|X|}^\leftarrow(1/4))] 
%\int_{(0,F_{|X|}^\leftarrow(1/4))}      \Bigl( 1 - F_{|X|}(0) \Bigr)\,d\leb(x)  \\*
      + \int_{[F_{|X|}^\leftarrow(1/4),\infty)} \Bigl( 1 - F_{|X|}(x) \Bigr)\,d\leb(x) 
\end{eqnarray*}
Since 
\begin{eqnarray*}
        \int_{[F_{|X|}^\leftarrow(1/4),\infty)} \Bigl( 1 - (D \circ F_{|X|})(x) \Bigr)\,d\leb(x)  
&  = &  \int_{[F_{|X|}^\leftarrow(1/4),\infty)} \Bigl( 1 -          F_{|X|} (x) \Bigr)\,d\leb(x) 
\end{eqnarray*}
we see that   $ \LL_Q^\Pichler = \LL^1 \neq \LL_{Q^E} = \LL_Q $. 
Consider finally a random variable   $ X $   satisfying 
\begin{eqnarray*}
        F_X(x)  
&  = &  \biggl( \frac{\beta}{-x} \biggr)^2 \chi_{(-\infty,-\beta)}(x) + \chi_{[-\beta,\infty)}(x) 
\end{eqnarray*}
for some   $ \beta\in(0,\infty) $. 
Then   $ -X $   has a Pareto distribution with finite expectation. 
This yields   $ X\in\LL^1=\LL_Q^\Pichler\subseteq\LL_Q $. 
Since   $ D(u) \geq (1/2)\,\sqrt{u}\,\chi_{[0,1/4)}(u) $,   
we obtain 
\begin{eqnarray*}
        \int_{(0,1)} |F_X^\leftarrow(u)|  \,dQ(u)  
&\geq&  \int_{(0,1)} (F_X^\leftarrow(u))^-\,dQ(u)  \\
&  = &  \int_{(-\infty,0)} (D \circ F_X)(x)\,d\leb(x)  \\
&\geq&  \int_{(-\infty,0)} \frac{1}{2}\,\sqrt{F_X(x)}\,\chi_{[0,1/4)}(F_X(x))\,d\leb(x)  \\
%&  = &  \int_{(-\infty,0)} \frac{1}{2}\,\sqrt{\biggl( \frac{\beta}{-x} \biggr)^2 \chi_{(-\infty,-\beta)}(x) + \chi_{[-\beta,\infty)}(x)}
%        \,\chi_{(-\infty,-2\beta)}(x)\,d\leb(x)  \\
&  = &  \int_{(-\infty,0)} \frac{1}{2}\,\biggl( \frac{\beta}{-x}\,\chi_{(-\infty,-\beta)}(x) + \chi_{[-\beta,\infty)}(x) \biggr)\,\chi_{(-\infty,-2\beta)}(x)\,d\leb(x)  \\
&  = &  \int_{(-\infty,-2\beta)} \frac{\beta}{-2\,x}\,d\leb(x) \\
&  = &  \frac{\beta}{2} \int_{(2\beta,\infty)} \frac{1}{z}\,d\leb(z)
\end{eqnarray*}
and hence   $ X\notin\LL_Q^\Acerbi $. 
Therefore, 
any two of the three classes   $ \LL_Q $,   $ \LL_Q^\Acerbi $   and   $ \LL_Q^\Pichler $   are distinct. 
\end{onelist}
%The inequalities in (2) and (3) follow from Example 
%\ref{quantile-E-e}. 
\end{examples}

%%%%%%%%%%%%%%%%%%%%%%%%%%%%%%%%%%%%%%%%%%%%%%%%%%%%%%%%%%%%%%%%%%%%%%%%%%%%%%%%%%%%%%%%%%%%%%%%%%%

\section*{Acknowledgement}
\label{acknowledgment}

The first and the last author gratefully acknowledge substantial discussions on risk measures with 
D{\'e}sir{\'e} D{\"o}rner, 
Andr{\'e} Neumann and 
Sarah Santo.

%%%%%%%%%%%%%%%%%%%%%%%%%%%%%%%%%%%%%%%%%%%%%%%%%%%%%%%%%%%%%%%%%%%%%%%%%%%%%%%%%%%%%%%%%%%%%%%%%%%

\section*{References}

\small
Acerbi C [2002]: 
Spectral measures of risk -- A coherent representation of subjective risk aversion. 
Journal of Banking \& Finance 26, 1505--1518. 

\smallskip
Aliprantis CD, Burkinshaw O [1990]: 
Principles of Real Analysis. 
Second Edition. 
Boston: Academic Press. 

\smallskip
Embrechts P, Wang R [2015]: 
Seven proofs for the subadditivity of expected shortfall.  
Dependence Modelling 3, 126--140.

\smallskip
Gzyl H, Mayoral S [2008]: 
On a relationship between distorted and spectral risk measures. 
Revista de Econom{\'i}a Financiera 15, 8--21. 

\smallskip
Pichler A [2013]: 
The natural Banach space for version independent risk measures. 
Insurance Mathematics and Economics 53, 405--415. 

\smallskip
Wang S, Dhaene J [1998]: 
Comonotonicity, correlation order and premium principles. 
Insurance Mathematics and Economics 22, 235--242. 

\smallskip
Wirch J, Hardy MR [2002]: 
Distortion risk measures -- Coherence and stochastic dominance. 
Proceedings of the Sixth Conference on Insurance Mathematics and Economics, Lisbon, July 2002.

\bigskip
\vfill\hspace*{\fill}\today

\end{document}